%% file: EstimationOfDistributedDelays.tex
\documentclass[a4paper,fleqn]{cas-dc}

\input{./preamble/packages}


\input{./preamble/commands}

\begin{document}
	\let\WriteBookmarks\relax
	\def\floatpagepagefraction{1}
	\def\textpagefraction{.001}

	\shorttitle{Identification of distributed delays}

	\shortauthors{Ritschel and Wyller}

	\title[mode=title]{
		An algorithm for distributed time delay identification based on a mixed Erlang kernel approximation and the linear chain trick
	}

	\author[1]{Tobias K. S. Ritschel}\ead{tobk@dtu.dk}
	\affiliation[1]{
		organization={Department of Applied Mathematics and Computer Science, Technical University of Denmark, DK-2800 Kgs. Lyngby, Denmark}}

	\cormark[1]
	\cortext[1]{Corresponding author: T.~K.~S.~Ritschel. Tel. +45 4525 3315.}

	\author[2]{John Wyller}\ead{john.wyller@nmbu.no}
	\affiliation[2]{
		organization={Department of Mathematics, Faculty of Science and Technology, Norwegian University of Life Sciences, 1432 Ås, Norway}}

	\begin{keywords}
		Delay differential equations \sep
		distributed time delays \sep
		mixed Erlang distribution \sep
		linear chain trick \sep
		system identification \sep
		parameter estimation \sep
		dynamical optimization \sep
		single-shooting \sep
		sensitivity equations.
	\end{keywords}

	\begin{abstract}
		Time delays are ubiquitous in industry and nature, and they significantly affect both transient dynamics and stability properties. Consequently, it is often necessary to identify and account for the delays when, e.g., designing a model-based control strategy. However, identifying delays in differential equations is not straightforward and requires specialized methods.
		Therefore, we propose an algorithm for identifying distributed delays in delay differential equations (DDEs) that only involves simulation of ordinary differential equations (ODEs).
		Specifically, we 1)~approximate the kernel in the DDEs (also called the memory function) by the probability density function of a mixed Erlang distribution and 2)~use the linear chain trick (LCT) to transform the resulting DDEs into ODEs.
		Finally, the parameters in the kernel approximation are estimated as the solution to a dynamical least-squares problem, and we use a single-shooting approach to approximate this solution.
		We demonstrate the efficacy of the algorithm using numerical examples that involve the logistic equation and a point reactor kinetics model of a molten salt nuclear fission reactor.
	\end{abstract}

	\maketitle

    \input{./tex/introduction} 
    \input{./tex/system} 
    \input{./tex/approximation} 
    \input{./tex/linear_chain_trick} 
    \input{./tex/algorithm} 
    \input{./tex/numerical_example} 
    \input{./tex/conclusions} 


\input{./EstimationOfDistributedDelays.bbl}

    \appendix

    \input{./appendix/derivation} 
    \input{./appendix/numerical_simulation} 
\end{document}

%% file: preamble/packages.tex
\usepackage[numbers]{natbib}

%% file: preamble/commands.tex
\newcommand{\incr}{\,\mathrm{d}}
\newcommand{\pdiff}[2]{\frac{\partial{#1}}{\partial{#2}}}

\newcommand{\R}{\mathbb{R}}
\newcommand{\N}{\mathbb{N}}

\newtheorem{assum}{Assumption}
\newtheorem{rem}{Remark}
\newtheorem{thm}{Theorem}

%% file: tex/introduction.tex
\section{Introduction}
Time delays exist in many industrial, biological, ecological, and environmental processes~\cite{Kolmanovskii:Myshkis:1992}, and understanding their nature is instrumental to accurately describing the dynamics of a system. They can be caused by many different phenomena, e.g., transport processes (flow in a pipe or a river), mixing dynamics (diffusion), growth kinetics, and communication delays. As they can significantly affect the dynamics and stability properties of a process~\cite{Niculescu:Gu:2004}, they must be accounted for when designing systems for monitoring, prediction, control, and optimization~\cite{Richard:2003}. Furthermore, although it is possible to model the individual delay processes in detail, it often requires more knowledge than is readily available. Instead, the process can be described by delay differential equations (DDEs)~\cite{Smith:2011} and the delay can be identified using measurements.

Delays are typically identified by minimizing a least-squares criterion~\cite{Rihan:2021} or maximizing a likelihood criterion (other criteria have been considered as well~\cite{Kozlowski:Kowalczuk:2015, Lin:etal:2015}). The resulting dynamical optimization problem is often solved approximately using an iterative numerical~\cite{Nocedal:Wright:2006} or metaheuristic~\cite{Burke:Kendall:2014} optimization method. Both types of methods involve repeated numerical simulation of the DDEs. However, numerical simulation of DDEs requires specialized methods~\cite{Bellen:Zennaro:2003}.
The reason is that, in contrast to ordinary differential equations (ODEs), the right-hand side of DDEs is a function of the state variables in one or more time intervals (distributed delays) or at multiple discrete points in time (absolute delays)~\cite{Balachandran:etal:2009}. Distributed delays involve an integral of all historical values of the states weighted by a kernel (also called a memory function). There exists general purpose software for numerical simulation of DDEs with absolute delays~\cite{Shampine:Thompson:2001, Yan:etal:2021}, but the solution is not differentiable with respect to the delay~\cite{Baker:Paul:1997}. Furthermore, to our knowledge, no such general purpose software exists for DDEs with distributed delays (although the Phoenix modeling language does contain some capabilities for distributed delays~\cite{Krzyzanski:etal:2018}). In conclusion, identification of delays in general nonlinear DDEs is not straightforward and requires specialized methods.
However, DDEs with distributed delays can be transformed into ODEs if the kernel is given by the probability density function of an Erlang distribution~\cite{Ibe:2014}. This transformation is referred to as the \emph{linear chain trick} (LCT)~\cite{MacDonald:1978, ponosov2004thew, Smith:2011}, and it relies on Leibniz' integral rule~\cite{Protter:Morrey:1985}.
Additionally, for an appropriate choice of weights, the cumulative distribution function of an infinite mixed Erlang distribution converges pointwise to the integral of any sufficiently smooth kernel as the (common) rate parameter goes to infinity~\cite{Tijms:1995, Cossette:etal:2016}.
This motivates the approximation of kernels by the probability density function of a finite mixed Erlang distribution with a finite common rate parameter.

Many authors have proposed methods for identifying time delays in systems with a specific structure, e.g., ARMAX models~\cite{Bedoui:Abderrahim:2015, Bedoui:Abderrahim:2022}, transfer function models~\cite{Wang:Zhang:2001}, and linear discrete-~\cite{Gu:etal:2021, Gu:Li:2023} and continuous-time state space models~\cite{Reiss:2005, Drakunov:2006, Belkoura:etal:2009}. Bj{\"{o}}rklund~\cite{Bjorklund:2003} presents a comprehensive survey of methods for linear systems. For differentially flat systems~\cite{Levine:2009}, the time delay can be estimated by transforming the DDEs into a set of algebraic equations~\cite{Schenkendorf:etal:2012, Schenkendorf:Mangold:2014}. Graphical methods for visually identifying absolute delays have also been proposed~\cite{Fowler:Kember:1993, Bunner:etal:1996, Ellner:etal:1997}.
For general nonlinear continuous-time systems with absolute delays, the DDEs can be discretized in order to apply methods for discrete-time systems~\cite{Hartung:Turi:2005, Kozlowski:Kowalczuk:2015}. Alternatively, the delayed states can be linearized with respect to the delay~\cite{Gao:etal:2022}, which also mitigates the issue of non-differentiability. Wang and Cao~\cite{Wang:Cao:2011} approximate the state variables using cubic B-splines and penalize the violation of the DDEs in the objective function. As non-differentiability is primarily a concern when using gradient-based numerical optimization methods, Tang and Guan~\cite{Tang:Guan:2009} and Jamilla et al.~\cite{Jamilla:etal:2021} use metaheuristic methods to approximate the solution to the involved dynamical optimization problem. Hauber et al.~\cite{Hauber:etal:2020} use the LCT to transform a set of DDEs into ODEs in order to estimate an absolute delay.

For systems with distributed delays, it is common to approximate either the kernel itself or the integral involving the historical states and the kernel. Torkamani et al.~\cite{Torkamani:etal:2013} use Gaussian quadrature to discretize the integral, and they use a spectral element method to approximate the solution to the DDEs. Krzyzanski et al.~\cite{Krzyzanski:etal:2018} identify distributed delays using the Phoenix software, which approximates the integral numerically. Zhou~\cite{Zhou:2016} proposes to either discretize the integral or approximate the kernel using splines. Murphy~\cite{Murphy:1990} uses linear splines to transform nonlinear deterministic DDEs into ODEs in order to identify both absolute and distributed delays.
Finally, Krzyzanski~\cite{Krzyzanski:2019} uses a truncated binomial expansion to approximate the integral for a kernel given by the probability density function of a gamma distribution. However, the LCT has not previously been used in the identification of distributed delays.

The main contribution of this work is an algorithm for approximately identifying distributed delays in DDEs, which only requires off-the-shelf software for numerical optimization and simulation of ODEs. First, we approximate the kernel using a finite mixed Erlang distribution with a finite common rate parameter. Next, we use the LCT to transform the DDEs into ODEs, and (for simplicity) we use a least-squares criterion to identify the parameters in the approximate kernel as well as model parameters and the initial states. We use a single-shooting approach~\cite{Binder:etal:2001} to approximate the solution to the least-squares dynamical optimization problem, and we implement it using \textsc{Matlab}'s~\cite{Mathworks:2023} \texttt{ode15s} for numerical simulation and \texttt{fmincon} for numerical optimization. The implementation is publicly available~\cite{Ritschel:Wyller:2023}.
Finally, we demonstrate the efficacy of the approach using numerical examples involving 1)~the logistic equation and 2)~a point reactor kinetics model of a molten salt nuclear fission reactor.

The remainder of the paper is organized as follows. In Section~\ref{Problem statement}, we describe the DDEs with distributed delays considered in this work, and in Section~\ref{Approximation}, we present the approximation based on the mixed Erlang distribution. Next, in Section~\ref{Linear chain trick}, we describe the LCT, and in Section~\ref{Estimation}, we present the proposed algorithm. The numerical examples are presented in Section~\ref{Examples}, and conclusions are given in Section~\ref{Conclusions}.

%% file: tex/system.tex
\section{System\label{Problem statement}}
In this work, we consider continuous-discrete systems in the form
\begin{subequations}\label{eq:system}
    \begin{align}
        \label{eq:system:x0}
        x(t) &= x_0(t), & t &\in (-\infty, t_0], \\
        \label{eq:system:x}
        \dot x(t) &= f(x(t), z(t), p), & t &\in [t_0, t_f], \\
        \label{eq:system:y}
        y(t_k) &= g(x(t_k), p), & k &= 0, \ldots, N,
    \end{align}
\end{subequations}
where the contribution from the historical values of the states  is given by
\begin{subequations}\label{eq:system:delay}
    \begin{align}
        \label{eq:system:z}
        z(t) &= \int\limits_{-\infty}^t \alpha(t - s) r(s)\,\mathrm ds, \\
        \label{eq:system:r}
        r(t) &= h(x(t), p).
    \end{align}
\end{subequations}
Here, $t \in \R$ is time, $t_0, t_f \in \R$ are the initial and final time, respectively, $x:\R \rightarrow \R^{n_x}$ are the states, $x_0:\R \rightarrow \R^{n_x}$ is the initial state function, $z:\R \rightarrow \R^{n_z}$ is the contribution from the delayed quantity $r:\R \rightarrow \R^{n_z}$, $\alpha:[0, \infty) \rightarrow [0, \infty)$ is the kernel, $y:\R \rightarrow \R^{n_y}$ is a vector of measurements, and $p\in\R^{n_p}$ are model parameters. Furthermore, $f:\R^{n_x}\times \R^{n_z} \times \R^{n_p}\rightarrow \R^{n_x}$ is the right-hand side function of the dynamical system, $g:\R^{n_x} \times \R^{n_p} \rightarrow \R^{n_y}$ is the measurement function, and $h:\R^{n_x} \times \R^{n_p} \rightarrow \R^{n_z}$ is the delay function.
Measurements are obtained from the system at $N+1 \in \N$ discrete points in time, $t_k \in \R$ for $k = 0, \ldots, N$, where $t_N = t_f$.
\begin{assum}
    The functions $f$ and $h$ are smooth on their respective domains such that the initial value problem~\eqref{eq:system}--\eqref{eq:system:delay} is well-posed. We refer to the book by Cushing~\cite{cushing2013} and the paper by Posonov et al.~\cite{ponosov2004thew} for more details.
\end{assum}
\begin{assum}\label{ass:alpha:cont:norm}
    The kernel $\alpha$ in~\eqref{eq:system:z} is normalized, i.e.,
    \begin{align}\label{eq:kernel:normalization}
        \int\limits_0^\infty \alpha(t)\,\mathrm dt &= 1.
    \end{align}
\end{assum}
%


%% file: tex/approximation.tex
\section{Mixed Erlang approximation\label{Approximation}}
We approximate the kernel $\alpha$ in~\eqref{eq:system:z} by $\hat \alpha^{(M)}: [0, \infty) \rightarrow [0, \infty)$ given by
\begin{align}\label{eq:kernel:approximation}
    \hat \alpha^{(M)}(t)
    &= \sum_{m=0}^M c_m \hat \alpha_m(t),
\end{align}
where $\hat \alpha_m: [0, \infty) \rightarrow [0, \infty)$ is the probability density function of an $m+1$'th order Erlang distribution~\cite{Ibe:2014}:
\begin{align}\label{eq:erlang:pdf}
    \hat \alpha_m(t) &= b_m t^m e^{-at}.
\end{align}
The normalization factor $b_m$ is
\begin{align}
    b_m &= \frac{a^{m+1}}{m!},
\end{align}
and the parameters in the approximation are the coefficients $c_m \in [0, 1]$ for $m = 0, \ldots, M$, the common rate parameter $a \in (0, \infty)$, and $M \in \N$. Furthermore, the coefficients must satisfy
\begin{align}\label{eq:delay:parameter:conditions}
    \sum_{m=0}^M c_m &= 1,
\end{align}
in order to ensure that
\begin{align}
    \int\limits_0^\infty \hat \alpha^{(M)}(t)\, \mathrm dt &= 1.
\end{align}
\begin{rem}
    The cumulative density function of an $m+1$'th order Erlang distribution is given by
    \begin{align}\label{eq:erlang:cdf}
        \hat \beta_m(t) &= \int\limits_0^t \hat \alpha_m(s)\,\mathrm ds = 1 - \frac{1}{a}\sum_{n=0}^m \hat \alpha_n(t).
    \end{align}
\end{rem}

\input{./tex/convergence}

%% file: tex/convergence.tex
\subsection{A note on convergence}
The following theorem states that the weights $\{c_m\}_{m=0}^M$ can be chosen such that the integral of the approximate kernel converges to that of the true kernel as the rate parameter $a$ and $M$ tend to infinity. The theorem follows directly from Theorem~2.9.1 in the book by Tijms~\cite{Tijms:1995} (see also Theorem~2 in~\cite{Cossette:etal:2016}).
\begin{thm}\label{thm:tijms}
    Let $\alpha:[0, \infty) \rightarrow [0, \infty)$ be a kernel that satisfies Assumption~\ref{ass:alpha:cont:norm} and define $\beta:[0, \infty) \rightarrow [0, 1]$ as
    \begin{align}
    	\beta(t) &= \int_0^t \alpha(s) \incr s.
    \end{align}
    Next, define $\hat \beta^{(M)}: [0, \infty) \rightarrow [0, 1]$ as
    \begin{align}
        \hat \beta^{(M)}(t)
        &= \int\limits_0^t \hat \alpha^{(M)}(s)\,\mathrm ds = \sum_{m=0}^M c_m \hat \beta_m(t),
    \end{align}
    where $c_m = \beta(t_{m+1}) - \beta(t_m)$, $a = 1/\Delta t$, and $t_m = m \Delta t$. Note that $\hat \beta_m$ depends on the rate parameter $a$. Then, $\hat \beta^{(\infty)}$ converges pointwise to $\beta$, i.e.,
    \begin{align}
        \lim_{\Delta t \rightarrow 0} \hat \beta^{(\infty)}(t) = \beta(t)
    \end{align}
    for any $t \in [0, \infty)$ where $\beta$ is continuous.
\end{thm}
\begin{rem}
    Scheff{\'{e}}~\cite{Scheffe:1947} proved that convergence of probability density functions (e.g,. $\hat \alpha^{(\infty)}$) implies convergence of cumulative distribution functions (e.g., $\hat \beta^{(\infty)}$). Later, Boos~\cite{Boos:1985} showed that the converse is also true if the probability density function is bounded and equicontinuous (see also the paper by Sweeting~\cite{Sweeting:1986}). The same holds if the probability density function converges uniformly~\cite{Rudin:1976}. However, it is outside the scope of this work to prove that $\hat \alpha^{(\infty)}$ converges to $\alpha$. Instead, we use the algorithm presented in Section~\ref{Estimation} to identify the coefficients $\{c_m\}_{m=0}^M$ in $\hat \alpha^{(M)}$ based on the measurements $\{y(t_k)\}_{k=0}^N$.
\end{rem}

%% file: tex/linear_chain_trick.tex
\section{The linear chain trick\label{Linear chain trick}}
We use the kernel approximation presented in Section~\ref{Approximation} to obtain the following system of DDEs:
\begin{subequations}\label{eq:lct:approximate:system:DDE}
    \begin{align}
    \dot{\hat x}(t) &= f(\hat x(t), \hat z(t), p), \\
        \hat z(t) &= \int\limits_{-\infty}^t \hat \alpha^{(M)}(t - s) \hat r(s)\, \mathrm ds, \\
        \hat r(t) &= h(\hat x(t), p).
    \end{align}
\end{subequations}
This system is in the same form as the original system~\eqref{eq:system:x} and~\eqref{eq:system:delay}, and $\hat x: \R \rightarrow \R^{n_x}$, $\hat z: \R \rightarrow \R^{n_z}$, and $\hat r: \R \rightarrow \R^{n_z}$ are approximations of $x$, $z$, and $r$, respectively.
We transform this system to a set of ODEs using the linear chain trick~\cite{MacDonald:1978}, as described in Appendix~\ref{app:derivation}. The resulting approximate system is
\begin{subequations}\label{eq:lct:approximate:system:ODE}
    \begin{align}
        \label{eq:system:approximate:x}
        \dot{\hat x}(t) &= f(\hat x(t), \hat z(t), p), \\
        \label{eq:system:approximate:Z}
        \dot{\hat Z}(t) &= A \hat Z(t) + B \hat r(t), \\
        \label{eq:system:approximate:z}
        \hat z(t) &= C \hat Z(t), \\
        \hat r(t) &= h(\hat x(t), p),
    \end{align}
\end{subequations}
where $\hat Z: \R \rightarrow \R^{(M+1) n_z}$ is a set of auxiliary state variables.
Furthermore, the system matrices $A\in \R^{(M+1) n_z \times (M+1) n_z}$, $B \in \R^{(M+1) n_z \times n_z}$, and $C \in \R^{n_z \times (M+1) n_z}$ in~\eqref{eq:system:approximate:Z}--\eqref{eq:system:approximate:z} are given by
\begin{subequations}\label{eq:system:approximate:matrices}
    \begin{align}
		A &= a
        \begin{bmatrix}
        -I \\
         I & -I \\
           &  \ddots & \ddots \\
           &         &      I & -I
        \end{bmatrix}, &
		B &= a
        \begin{bmatrix}
            I \\
            \phantom{1} \\
            \phantom{\vdots} \\
            \phantom{1}
        \end{bmatrix}, \\
		C &=
        \begin{bmatrix}
            c_0 I & c_1 I & \cdots & c_M I
        \end{bmatrix},
    \end{align}
\end{subequations}
where $I \in \R^{n_z \times n_z}$ is an identity matrix.

%% file: tex/algorithm.tex
\section{Algorithm\label{Estimation}}
We formulate the delay identification problem as a constrained dynamical parameter estimation problem with a least-squares objective function. The dynamical constraints are the approximate system of ODEs~\eqref{eq:lct:approximate:system:ODE} and, for completeness, we also estimate model parameters and the initial state.
Furthermore, we use a single-shooting approach~\cite{Binder:etal:2001} to approximate the solution to the estimation problem, we compute the involved first-order sensitivities using a forward approach, and we briefly describe a \textsc{Matlab}~\cite{Mathworks:2023} implementation of the algorithm.
\begin{rem}
    For simplicity, we use a least-squares criterion in the estimation problem. However, other criteria, e.g., maximum likelihood criteria, are also applicable.
\end{rem}
\begin{rem}
    The solution to the dynamical parameter estimation problem can also be approximated using, e.g., multiple-shooting or simultaneous collocation methods~\cite{Binder:etal:2001}.
\end{rem}
\input{./tex/dynamical_least_squares}
\input{./tex/single_shooting}
\input{./tex/sensitivities}
\input{./tex/implementation}

%% file: tex/dynamical_least_squares.tex
\subsection{Dynamical least-squares problem}
For a given value of the kernel parameter $M$, the estimates of 1)~the model parameters, $\hat p \in \R^{n_p}$, 2)~the $n_q = M+2 \in \N$ remaining kernel parameters, $\hat q \in \R^{n_q}$, i.e., $\{c_m\}_{m=0}^M$ and $a$, and 3)~the initial states, $\hat x_0 \in \R^{n_x}$, are the solution to the dynamical least-squares optimization problem
\begin{align}\label{eq:least:squares:obj}
    \min_{\hat p, \hat q, \hat x_0} \quad \phi\left(\{y(t_k)\}_{k=0}^N, \{\hat y(t_k)\}_{k=0}^N\right),
\end{align}
subject to
\begin{subequations}\label{eq:least:squares}
    \begin{align}
        \label{eq:least:squares:x0}
        \hat x(t_0) &= \hat x_0, \\
        \label{eq:least:squares:Z0}
        \hat Z(t_0) &= \hat Z_0(\hat x_0, \hat p), \\
        \label{eq:least:squares:x}
        \dot{\hat x}(t) &= f(\hat x(t), \hat z(t), \hat p), & t &\in [t_0, t_f], \\
        \label{eq:least:squares:Z}
        \dot{\hat Z}(t) &= A(\hat q) \hat Z(t) + B(\hat q) \hat r(t), & t &\in [t_0, t_f], \\
        \label{eq:least:squares:z}
        \hat z(t) &= C(\hat q) \hat Z(t), & t &\in [t_0, t_f], \\
        \label{eq:least:squares:r}
        \hat r(t) &= h(\hat x(t), \hat p), & t &\in [t_0, t_f], \\
        \label{eq:least:squares:y}
        \hat y(t_k) &= g(\hat x(t_k), \hat p), & k &= 0, \ldots, N, \\
        \label{eq:least:squares:c}
        w^T \hat q &= 1, \\
        \label{eq:least:squares:p}
        p_{\min} &\leq \hat p \leq p_{\max}, \\
        \label{eq:least:squares:q}
        q_{\min} &\leq \hat q \leq q_{\max}, \\
        \label{eq:least:squares:x:bounds}
        x_{\min} &\leq \hat x_0 \leq x_{\max}.
    \end{align}
\end{subequations}
The objective function $\phi: \R^{n_y} \times \dots \times \R^{n_y} \rightarrow \R$ penalizes the squared deviations from the measurements:
\begin{multline}
    \phi(\{y(t_k)\}_{k=0}^N, \{\hat y(t_k)\}_{k=0}^N) = \\
    \frac{1}{2} \sum_{k=0}^N (y(t_k) - \hat y(t_k))^T (y(t_k) - \hat y(t_k)).
\end{multline}
We assume that the system is in steady state up until time $t_0$, and~\eqref{eq:least:squares:x0}--\eqref{eq:least:squares:Z0} are initial conditions. The initial value $\hat Z_0: \R^{n_x} \times \R^{n_p} \rightarrow \R^{(M+1) n_z}$ is a function of the estimated initial state $\hat x_0$ and the parameter estimate $\hat p$ (see Remark~\ref{rem:steady:state} in Appendix~\ref{app:derivation}). Furthermore, the dynamical constraints \eqref{eq:least:squares:x}--\eqref{eq:least:squares:r} are the transformed approximate system of ODEs presented in Section~\ref{Linear chain trick}, \eqref{eq:least:squares:y} is the measurement equation, \eqref{eq:least:squares:c} represents the condition~\eqref{eq:delay:parameter:conditions} on $\{c_m\}_{m=0}^M$, and \eqref{eq:least:squares:p}--\eqref{eq:least:squares:x:bounds} are bound constraints with given values of $p_{\min}, p_{\max} \in \R^{n_p}$, $q_{\min}, q_{\max} \in \R^{n_q}$, and $x_{\min}, x_{\max} \in \R^{n_x}$. The actual measurements, $\{y(t_k)\}_{k=0}^N$, are also given. The element of $w \in \R^{n_q}$ corresponding to $a$ is 0 and all other elements are 1. Furthermore, we ensure that $a$ is positive using an arbitrary positive lower bound. Finally, the functions $A: \R^{n_q} \rightarrow \R^{(M+1) n_z \times (M+1) n_z}$, $B: \R^{n_q} \rightarrow \R^{(M+1) n_z \times n_z}$, and $C: \R^{n_q} \rightarrow \R^{n_z \times (M+1) n_z}$ are given by~\eqref{eq:system:approximate:matrices}.
\begin{rem}
    If the system cannot be assumed to be in steady state up until time $t_0$, the initial state function can be approximated by a parametrized function, the integral in~\eqref{eq:identification:initial:conditions:zm0} defining the elements of $\hat Z_0$ can be approximated using quadrature, and the parameters in the parametrized function can be estimated instead of $\hat x_0$.
\end{rem}

%% file: tex/single_shooting.tex
\subsection{Single-shooting}\label{app:single:shooting}
In the single-shooting approach, we transform the dynamical optimization problem~\eqref{eq:least:squares:obj}--\eqref{eq:least:squares} to the nonlinear program (NLP)
\begin{subequations}\label{eq:single:shooting}
    \begin{align}
        \min_{\hat p, \hat q, \hat x_0} \hspace{22pt} &\psi(\hat p, \hat q, \hat x_0), \\
        \text{subject to} \hspace{10pt} & \text{\eqref{eq:least:squares:c}--\eqref{eq:least:squares:x:bounds}},
    \end{align}
\end{subequations}
where the objective function, $\psi: \R^{n_p} \times \R^{n_q} \times \R^{n_x} \rightarrow \R$, is given by
\begin{multline}\label{eq:single:shooting:obj}
    \psi(\hat p, \hat q, \hat x_0) = \\ \Big\{\phi\left(\{y(t_k)\}_{k=0}^N, \{\hat y(t_k)\}_{k=0}^N\right) : \text{\eqref{eq:least:squares:x0}--\eqref{eq:least:squares:y}}\Big\}.
\end{multline}
That is, given $\hat p$, $\hat q$, and $\hat x_0$, $\psi$ is the objective function $\phi$ evaluated using the solution to the initial value problem~\eqref{eq:least:squares:x0}--\eqref{eq:least:squares:r} and the measurement equation~\eqref{eq:least:squares:y}.

%% file: tex/sensitivities.tex
\subsection{Sensitivities}\label{app:sensitivities}
For single-shooting approaches, the objective function $\psi$ can be a highly nonlinear function of the decision variables. Consequently, it can be difficult to accurately approximate the derivatives of $\psi$ using, e.g., finite difference methods~\cite{LeVeque:2007}. Therefore, and as many numerical optimization methods are gradient-based~\cite{Nocedal:Wright:2006}, we use a forward sensitivity-based approach~\cite[Chap.~I.14]{Hairer:etal:1993} to compute the first-order derivatives.

First, we introduce $\theta \in \R^{n_p + n_q + n_x}$ given by
\begin{align}
    \theta &=
    \begin{bmatrix}
        \hat p \\ \hat q \\ \hat x_0
    \end{bmatrix}.
\end{align}
Then, the derivatives of the objective function~\eqref{eq:single:shooting:obj} in the NLP are
\begin{align}\label{eq:sensitivities:objective:function}
    \pdiff{\psi}{\theta_i}(\theta)
    &= -\sum_{k=0}^N (y(t_k) - \hat y(t_k))^T S_{y, i}(t_k).
\end{align}
Furthermore, $S_{y, i}: \R \rightarrow \R^{n_y}$ are the sensitivities (i.e., derivatives) of $\hat y$ with respect to the $i$'th element of $\theta$:
\begin{align}
    S_{y, i} &= \pdiff{\hat y}{\theta_i}.
\end{align}
We also introduce the sensitivities $S_{x, i}: \R \rightarrow \R^{n_x}$, $S_{Z, i}: \R \rightarrow \R^{(M+1) n_z}$, and $S_{z, i}, S_{r, i}: \R \rightarrow \R^{n_z}$ given by
\begin{align}
    S_{x, i} &= \pdiff{\hat x}{\theta_i}, &
    S_{Z, i} &= \pdiff{\hat Z}{\theta_i}, &
    S_{z, i} &= \pdiff{\hat z}{\theta_i}, &
    S_{r, i} &= \pdiff{\hat r}{\theta_i}.
\end{align}
These sensitivities are obtained as the solution to the initial value problem
\begin{subequations}\label{eq:sensitivity:equations}
    \begin{align}
        S_{x, i}(t_0)
        &= \pdiff{\hat x_0}{\theta_i}, \\
        S_{Z, i}(t_0)
        &= \pdiff{\hat Z_0}{\theta_i}(\hat x_0, \hat p), \\
        \dot S_{x, i}(t)
        &= \pdiff{f}{x} S_{x, i}(t) + \pdiff{f}{z} S_{z, i}(t) + \pdiff{f}{\theta_i}, \\
        \dot S_{Z, i}(t)
        &= \pdiff{A}{\theta_i} \hat Z(t) + A(\hat q) S_{Z, i}(t) \nonumber \\
        &+ \pdiff{B}{\theta_i} \hat r(t) + B(\hat q) S_{r, i}(t), \\
        S_{z, i}(t)
        &= \pdiff{C}{\theta_i} \hat Z(t) + C(\hat q) S_{Z, i}(t), \\
        S_{r, i}(t)
        &= \pdiff{h}{x} S_{x, i}(t) + \pdiff{h}{\theta_i},
    \end{align}
\end{subequations}
for $t \in [t_0, t_f]$. Using the solution to this initial value problem, we compute the necessary sensitivities by
\begin{align}\label{eq:sensitivity:equations:y}
    S_{y, i}(t_k)
    &= \pdiff{g}{x} S_{x, i}(t_k) + \pdiff{g}{\theta_i}, & k &= 0, \ldots, N.
\end{align}
For brevity of notation, we have omitted the arguments of the Jacobians. In~\eqref{eq:sensitivity:equations}, the Jacobians of $f$ are evaluated in $\hat x(t)$, $\hat z(t)$, and $\hat p$ and the Jacobians of $h$ are evaluated in $\hat x(t)$ and $\hat p$. Furthermore, the Jacobians of $g$ in~\eqref{eq:sensitivity:equations:y} are evaluated in $\hat x(t_k)$ and $\hat p$. As $A$, $B$, and $C$ are linear in the elements of $\hat q$, their derivatives are constant. Finally, in~\eqref{eq:sensitivities:objective:function} and~\eqref{eq:sensitivity:equations}--\eqref{eq:sensitivity:equations:y}, $\hat x$, $\hat Z$, $\hat z$, and $\hat r$ are the solutions to the initial value problem~\eqref{eq:least:squares:x0}--\eqref{eq:least:squares:r} and $\hat y$ is given by~\eqref{eq:least:squares:y}.

%% file: tex/implementation.tex
\subsection{Implementation}
We implement the algorithm described in this section using \textsc{Matlab}. Specifically, we use \texttt{ode15s} to approximate the solution to the initial value problem combining~\eqref{eq:least:squares:x0}--\eqref{eq:least:squares:r} and~\eqref{eq:sensitivity:equations}, i.e., the initial value problems involving the approximate system of ODEs and the corresponding sensitivity equations. Furthermore, we provide the required Jacobian of the right-hand side of the differential equations. We derive the Jacobian analytically, and we implement it as a sparse matrix. We use \texttt{fmincon}'s interior point algorithm to approximate the solution to the NLP~\eqref{eq:single:shooting} and, as described in Section~\ref{app:sensitivities}, we provide the gradient of the objective function $\psi$.  However, we do not provide the Hessian matrix. Instead, it is approximated by \texttt{fmincon} using a BFGS approximation~\cite{Nocedal:Wright:2006}.

%% file: tex/numerical_example.tex
\section{Numerical examples\label{Examples}}
In this section, we present two numerical examples that demonstrate the efficacy of the algorithm described in Section~\ref{Estimation}. The first example involves the logistic equation with a time-varying carrying capacity and the second involves a point reactor kinetics model of a molten salt nuclear fission reactor.

\input{./tex/logistic_equation} 
\input{./tex/nuclear_fission}

%% file: tex/logistic_equation.tex
\subsection{The logistic equation}
We consider the logistic equation
\begin{align}
    \dot N(t) &= \kappa N(t) \left(1 - \frac{\tilde N(t)}{K(t)}\right),
\end{align}
where $\kappa \in [0, \infty)$ is the logistic growth rate, $N: \R \rightarrow [0, \infty)$ is the relative population density, and $\tilde N: \R \rightarrow [0, \infty)$ is the delayed population density given by
\begin{align}\label{eq:logistic:delay}
    \tilde N(t) &= \int\limits_{-\infty}^{t} \alpha(t - s) N(s)\,\mathrm ds.
\end{align}
Furthermore, $K: \R \rightarrow (0, \infty)$ is the time-varying carrying capacity, which is given by
\begin{align}
    K(t) &= \left(1 + A_1 \sin(2\pi \omega_1 t) + A_2 \sin(2\pi \omega_2 t)\right) \bar K,
\end{align}
where $A_1, A_2 \in [0, \infty)$ are amplitudes, $\omega_1, \omega_2 \in [0, \infty)$ are frequencies, and $\bar K \in (0, \infty)$ is the nominal carrying capacity.
Finally, the kernel is the probability density function of a mixed folded normal distribution, i.e.,
\begin{align}\label{eq:logistic:equation:kernel}
    \alpha(t) &= \gamma_1 F(t; \mu_1, \sigma_1) + \gamma_2 F(t; \mu_2, \sigma_2),
\end{align}
where $\gamma_1, \gamma_2 \in [0, 1]$ are weights, $\mu_1, \mu_2 \in \R$ are location parameters, $\sigma_1, \sigma_2 \in (0, \infty)$ are scale parameters, and $F: [0, \infty) \times \R \times (0, \infty) \rightarrow [0, \infty)$ is given by
\begin{align}
    F(t; \mu, \sigma) &= \frac{\exp\left(-\frac{1}{2}\left(\frac{t - \mu}{\sigma}\right)^2\right) + \exp\left(-\frac{1}{2}\left(\frac{t + \mu}{\sigma}\right)^2\right)}{\sqrt{2 \pi}\, \sigma}.
\end{align}
Table~\ref{tab:logistic:parameters} shows the parameter values. The frequencies represent seasonal and monthly variations, and we choose the kernel such that there are no large undamped oscillations. A complete stability and bifurcation analysis is outside the scope of this work. However, in previous work, such an analysis has been carried out for systems with distributed delays using the Routh-Hurwitz criterion~\cite{Nordbo2007}.
We use the numerical method described in Appendix~\ref{sec:numerical:simulation}, with 150 time steps per day, to obtain daily measurements of the population density over a 2-year period. The initial density is $N(t) = N_0 = 0.9$ for $t \in (-\infty, t_0]$. Furthermore, we use a memory horizon of $\Delta t_h = 24$~mo, and we use a tolerance of~$10^{-12}$ when solving the involved residual equations using \textsc{Matlab}'s \texttt{fsolve}.

The objective is to estimate the kernel, the initial population density, $N_0$, and the growth rate, $\kappa$. The lower bounds on $\{c_m\}_{m=0}^M$, $N_0$, and $\kappa$ are~0 (in the respective units) and the lower bound on $a$ is~0.5~mo$^{-1}$. The upper bound on $N_0$ is~10, the upper bound on $\kappa$ is~10~mo$^{-1}$, and there are no upper bounds on the remaining parameters. The initial guess of $c_m$ is $1/(M+1)$ for $m = 0, \ldots, M$, and the initial guesses of $a$, $N_0$, and $\kappa$ are 20~mo$^{-1}$, 0.7, and 3~mo$^{-1}$, respectively.
Furthermore, although it does not change the solution to the optimization problem, we scale the objective function by a factor of $10^6$ in order to improve the convergence of \texttt{fmincon}. Finally, we use an optimality tolerance of~$10^{-3}$ in \texttt{fmincon} and an absolute and relative tolerance of~$10^{-8}$ in \texttt{ode15s}.
\begin{table}[t]
    \tiny
    \centering
    \caption{Values of the parameters in the logistic equation.}
    \label{tab:logistic:parameters}
    \begin{tabular}{cccccc}
        \hline
        \multicolumn{6}{c}{Model parameters} \\
        \hline
        $\kappa$~[1/mo] & $A_1$~[--] & $A_2$~[--] & $\omega_1$~[1/mo] & $\omega_2$~[1/mo] & $\bar K$~[--]  \\
        4 & 0.01 & 0.005 & 1/12 & 1 & 1 \\
        \hline
        \multicolumn{6}{c}{Kernel parameters} \\
        \hline
        $\gamma_1$~[--] & $\gamma_2$~[--] & $\mu_1$~[mo] & $\mu_2$~[mo] & $\sigma_1$~[mo] & $\sigma_2$~[mo] \\
        0.5 & 0.5 & 0.35 & 0.45 & 0.06 & 0.12 \\
        \hline
    \end{tabular}
\end{table}

\begin{figure*}[t]
    \centering
    \includegraphics[width=0.46\textwidth]{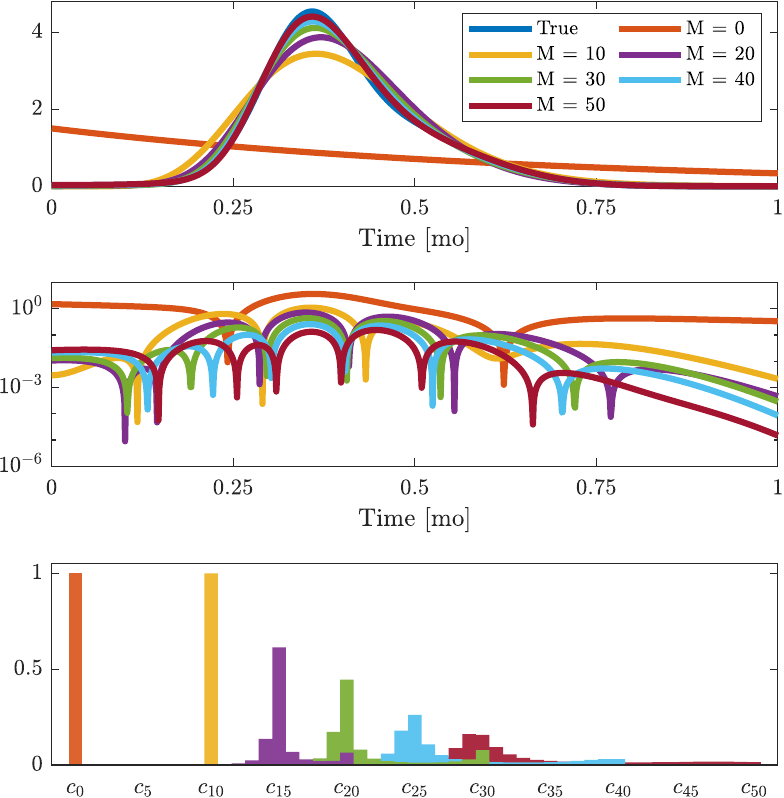}~
    \includegraphics[width=0.47\textwidth]{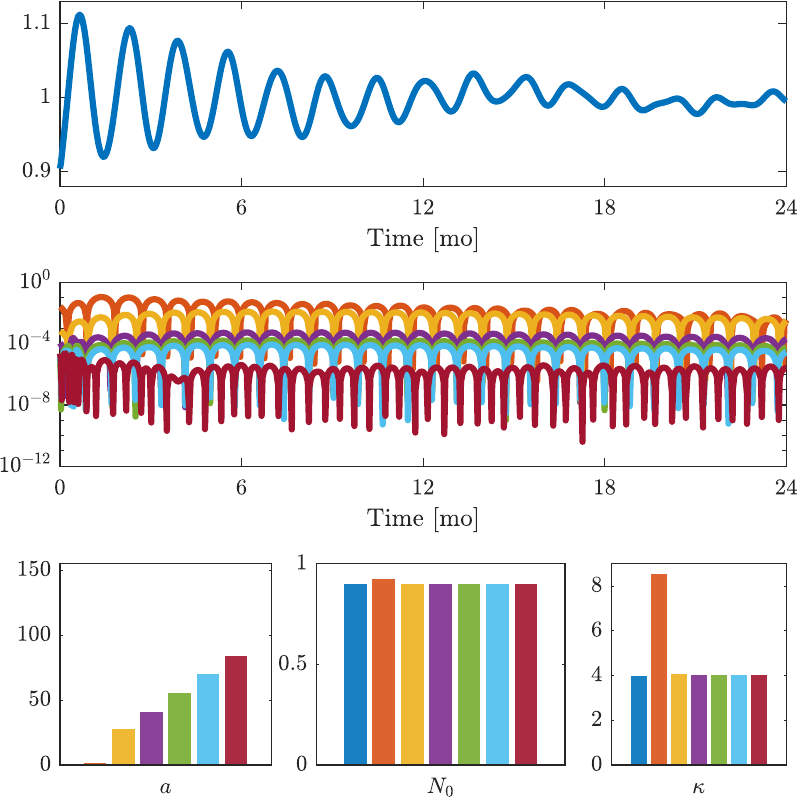}
    \caption{Estimation results for the logistic equation with a distributed delay. Left column: The true (hardly visible) and estimated kernels (top), the corresponding absolute error (middle), and the estimates of the coefficients $\{c_m\}_{m=0}^M$ (bottom). Right column: The population density for the true parameters (top), the absolute difference in population density for the true and estimated parameters (middle), and the true (when applicable) and estimated values of $a$, $N_0$, and $\kappa$ (bottom). The colors are consistent across the figure.}
    \label{fig:logistic:distributed}
\end{figure*}

Fig.~\ref{fig:logistic:distributed} shows the estimation results for $M = 0, 10, \ldots, 50$. The initial state and the growth rate are well estimated for the nonzero values of $M$. As $M$ increases, the estimate of $a$ increases as well, and the maximum kernel and population density errors decrease. For $M$ equal to~0 and~$10$, $c_M$ is estimated to be 1 and almost 1, respectively, whereas there are multiple nonzero coefficient estimates for the larger values of $M$. In summary, the results indicate that the algorithm can identify the kernel with high precision when $M$ is chosen sufficiently large.

\subsubsection{Absolute delay}
For completeness, we also test the algorithm using a system with an absolute delay, i.e., we replace~\eqref{eq:logistic:delay} by
\begin{align}
    \tilde N(t) &= N(t - \tau),
\end{align}
where the absolute delay is $\tau = 0.35$~mo. We simulate the system using \textsc{Matlab}'s \texttt{dde23} with absolute and relative tolerances equal to~$10^{-8}$.
In this case, we scale the objective function by $10^5$ and we use an optimality tolerance of $10^{-4}$ in \texttt{fmincon}. All other values remain unchanged.

\begin{figure*}[t]
    \centering
    \includegraphics[width=0.46\textwidth]{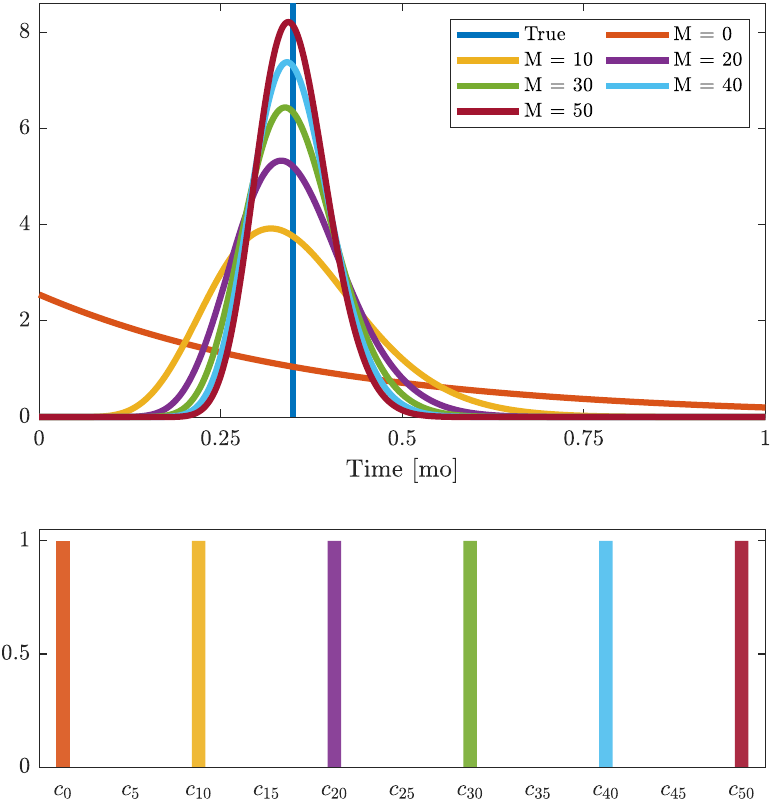}~
    \includegraphics[width=0.47\textwidth]{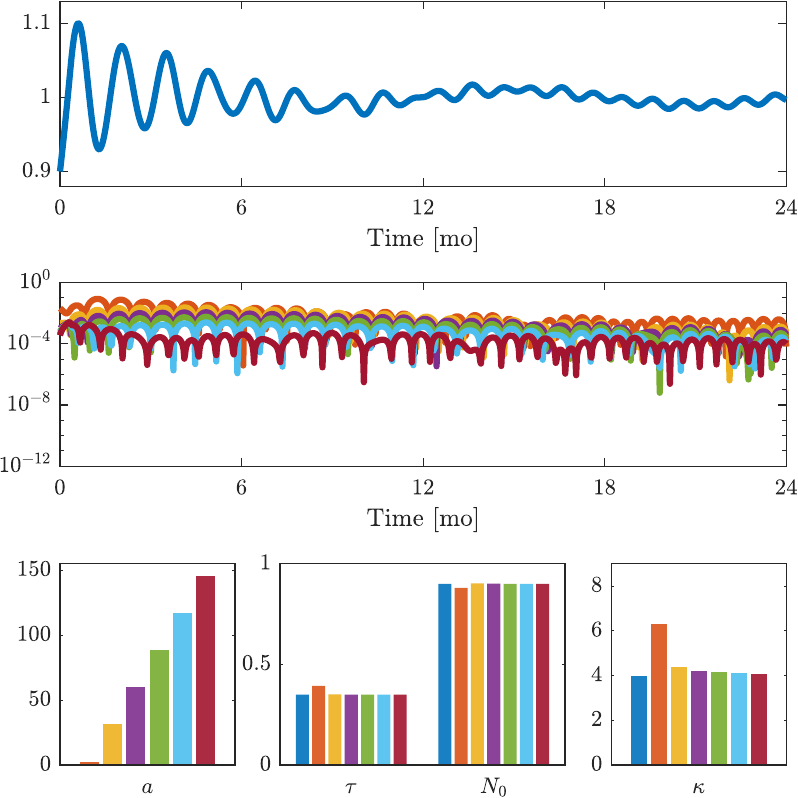}
    \caption{Estimation results for the logistic equation with an absolute delay. Left column: The true delay and the estimated kernels (top) and the estimates of the coefficients $\{c_m\}_{m=0}^M$ (bottom). Right column: The population density for the true parameters (top), the absolute difference in population density for the true and the estimated parameters (middle), and the true (when applicable) and estimated values of $a$, $\tau$, $N_0$, and $\kappa$ (bottom). The colors are consistent across the figure.}
    \label{fig:logistic:absolute}
\end{figure*}

We repeat the previous numerical experiments, and Fig.~\ref{fig:logistic:absolute} shows the results. As $M$ increases, the approximate kernel becomes more narrow around the true delay, and in all cases, $c_m$ is estimated to be almost~1 for $m = M$. Furthermore, the mean of the estimated kernel is an accurate estimate of the delay for the nonzero values of $M$, i.e., $\hat \tau = \frac{1}{a}\sum_{m=0}^M c_m (m+1) \approx \tau$. However, as expected, the difference in the population density is significantly larger than when estimating the distributed delay. Nevertheless, the results demonstrate that the algorithm can also be used to estimate absolute delays with high precision.

%% file: tex/nuclear_fission.tex
\subsection{Nuclear fission}
We consider a point reactor kinetics model of a molten salt nuclear reactor~\cite{Duderstadt:Hamilton:1976, Wooten:Powers:2018}, which accounts for 1)~neutron emission from fission events and the decay of fission products, 2)~circulation of the reactor fuel outside of the core, and 3)~temperature-dependent reactivity (negative feedback).
The fission reactions are in the form
\begin{align}
    {}_0^1{\mathrm n} + {}_{\phantom{2}92}^{235}{\mathrm U} \overset{\substack{\text{neutron}\\\text{absorp.}}}{\xrightarrow{\hspace{25pt}}} {}_{\phantom{2}92}^{236}{\mathrm U} \overset{\text{fission}}{\xrightarrow{\hspace{25pt}}} \sum_{i=1}^{N_p} \nu_i\,{}_{z_i}^{a_i}{\mathrm E_i} + \nu_0\,{}_0^1{\mathrm n}
\end{align}
where a neutron, ${}_0^1{\mathrm n}$, is absorbed by the uranium isotope ${}_{\phantom{2}92}^{235}{\mathrm U}$ to form the isotope ${}_{\phantom{2}92}^{236}{\mathrm U}$ which fissions into $N_p \in \N$ products. The $i$'th product is ${}_{z_i}^{a_i}{\mathrm E_i}$, where $a_i, z_i \in \N$ are the atomic mass and atomic number, respectively. Furthermore, $\nu_i \in \N$ for $i = 0, \ldots, N_p$ are stochiometric coefficients. A fraction of certain products undergoes negative beta decay followed by a neutron emission, i.e., the reaction is (irrelevant byproducts are omitted)
\begin{align}
    {}_{z_i}^{a_i}{\mathrm E_i} \overset{\substack{\text{neg. beta}\\\text{decay}}}{\xrightarrow{\hspace{25pt}}} {}_{z_i+1}^{\phantom{-1}a_i}{\mathrm E_i^*} \overset{\substack{\text{neutron}\\\text{emission}}}{\xrightarrow{\hspace{25pt}}} {}_{z_i+1}^{a_i-1}{\mathrm E_i^*} + {}_0^1{\mathrm n},
\end{align}
where ${}_{z_i+1}^{\phantom{-1}a_i}{\mathrm E_i^*}$ is a different element than ${}_{z_i}^{a_i}{\mathrm E_i}$ because the atomic number is one higher. The neutron-emitting products are called \emph{delayed neutron precursors} and in the model, they are grouped into $N_g = 6$ groups based on their half-lives.
The model consists of $n = N_g + 1$ mass balance equations as well as a reactivity equation derived from an energy balance. The derivation assumes that 1)~the power production is proportional to the neutron concentration, 2)~the time derivative of the reactivity is negatively proportional to that of the temperature, 3)~the heat loss due to external circulation of the fuel is insignificant, and 4)~the heat capacity of the reactor core is constant in time.
We formulate the model in a general form used to describe stirred tank reactors~\cite{Wahlgreen:etal:2020}.
It describes the temporal evolution of the reactivity $\rho: \R \rightarrow \R$ and the concentrations $C_i: \R \rightarrow [0, \infty)$ of the neutron-emitting fraction of the delayed precursor groups ($i = 1, \ldots, N_g$) and the neutrons ($i = n$):
\begin{subequations}
    \begin{align}
        \dot C_i(t) &= (C_{i, in}(t) - C_i(t)) D + R_i(t), \\
        \dot C_n(t) &= R_n(t), \\
        \dot \rho(t) &= -\kappa H C_n(t),
    \end{align}
\end{subequations}
for $i = 1, \ldots, N_g$. The dilution rate $D \in [0, \infty)$ is the ratio between the volumetric inlet and outlet flow rate and the reactor core volume, which is equal to the inverse of the fuel residence time $\tau_c \in (0, \infty)$. Furthermore, $\kappa \in [0, \infty)$ is the reactivity proportionality constant and $H \in [0, \infty)$ is the ratio between the power production proportionality constant and the heat capacity of the reactor core. The production rate $R: \R \rightarrow \R^n$ is given by
\begin{align}
    R(t) &= S^T r(t),
\end{align}
where $S \in \R^{N_r \times n}$ contains the stochiometric coefficients of the $N_r = n$ reactions and $r: \R \rightarrow [0, \infty)^{N_r}$ contains reaction rates:
\begin{align}
    S &=
    \begin{bmatrix}
        -1      &         &             & 1               \\
                &  \ddots &             & \vdots          \\
                &         & -1          & 1               \\
        \beta_1 & \cdots  & \beta_{N_g} & \rho(t) - \beta
    \end{bmatrix}, &
    r(t) &=
    \begin{bmatrix}
        \lambda_1 C_1(t) \\ \vdots \\ \lambda_{N_g} C_{N_g}(t) \\ C_n(t)/\Lambda
    \end{bmatrix}.
\end{align}
Here, $\lambda_i, \beta_i, \Lambda \in [0, \infty)$ for $i = 1, \ldots, N_g$ are decay constants, delayed neutron fractions, and the mean neutron generation time, respectively, and $\beta \in [0, \infty)$ is the sum of $\beta_i$ for $i = 1, \ldots, N_g$.
Next, the inlet concentration $C_{i, in}: \R \rightarrow [0, \infty)$ is
\begin{align}
    C_{i, in}(t) &= \delta_i \tilde C_i(t), & i &= 1, \ldots, N_g,
\end{align}
where $\tilde C_i: \R \rightarrow [0, \infty)$ is given by
\begin{align}\label{eq:nuclear:fission:delay}
    \tilde C_i(t) &= \int\limits_{-\infty}^t \alpha(t - s) C_i(s)\,\mathrm ds, & i &= 1, \ldots, N_g.
\end{align}
For simplicity, we assume that the decay factor $\delta_i = e^{-\lambda_i \tau_\ell} \in (0, \infty)$ only depends on the \emph{average} time spent outside of the reactor core, $\tau_\ell \in [0, \infty)$. Specifically, $\tau_\ell$ is the mean associated with the kernel $\alpha$, and $\alpha$ is given by~\eqref{eq:logistic:equation:kernel}. In this example, the distributed delay is chosen arbitrarily. However, it can represent that, e.g., friction between the molten salt and the pipe wall in the external loop causes the velocity profile to vary across the cross section of the pipe. In contrast, an absolute delay would represent plug flow. Table~\ref{tab:nuclear:fission:parameters} shows the parameter values, which (except for $\tau_c$ and $\tau_\ell$) are taken from~\cite{Leite:etal:2016}.

The objective is to estimate the kernel, the initial concentrations $C_{i, 0} = 1$~cm$^{-3}$ for $i = 1, \ldots, n$, the initial reactivity $\rho_0 = 1.5 \beta$, and the reactivity proportionality constant, $\kappa$. We use the method described in Appendix~\ref{sec:numerical:simulation} with 1,000 time steps per second to obtain 100~measurements per second of the concentrations over a period of 25~s. Furthermore, we use a memory horizon of $\Delta t_h = 25$~s and a tolerance of $10^{-12}$ when solving the involved residual equations. As the concentrations span several orders of magnitude, we use logarithmic concentrations in the measurement equation.
We scale the objective function by a factor of $10$, we use an optimality tolerance of $10^{-5}$ in \texttt{fmincon}, and we use an absolute and relative tolerance of $10^{-8}$ in \texttt{ode15s}.
The lower bounds on $\{c_m\}_{m=0}^M$, $\{C_{i, 0}\}_{i=1}^n$, and $\kappa$ are 0 (in the respective units) and the lower bound on $a$ is 7.5~s$^{-1}$. The upper bound on $\kappa$ is $10^{-4}$~K$^{-1}$, and there are no upper bounds on the remaining parameters. The initial guess of $c_m$ is 1 for $m = M$ and $10^{-8}$ for $m = 0, \ldots, M-1$. The initial guesses of $a$ and $\kappa$ are 25~s$^{-1}$ and $4 \cdot 10^{-5}$~K$^{-1}$, respectively, and the initial guesses of the initial concentrations and reactivity are 10~cm$^{-3}$ and $\beta$, respectively.

\begin{figure*}[t]
    \centering
    \includegraphics[width=0.455\textwidth]{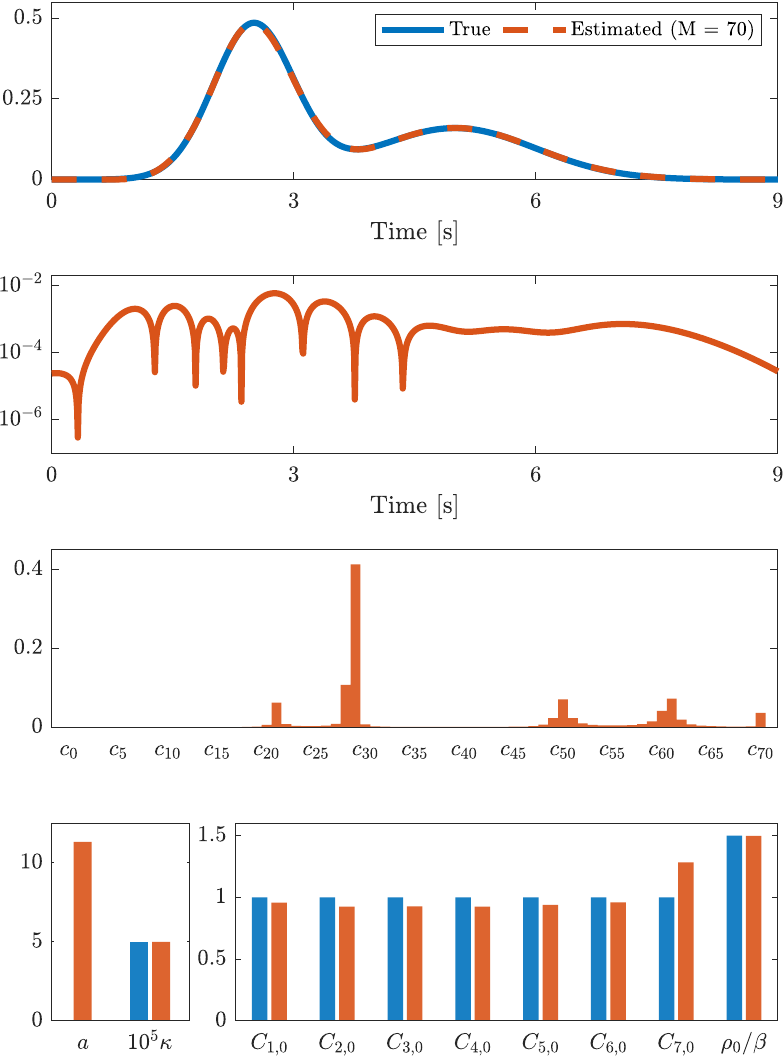}~
    \includegraphics[width=0.49\textwidth]{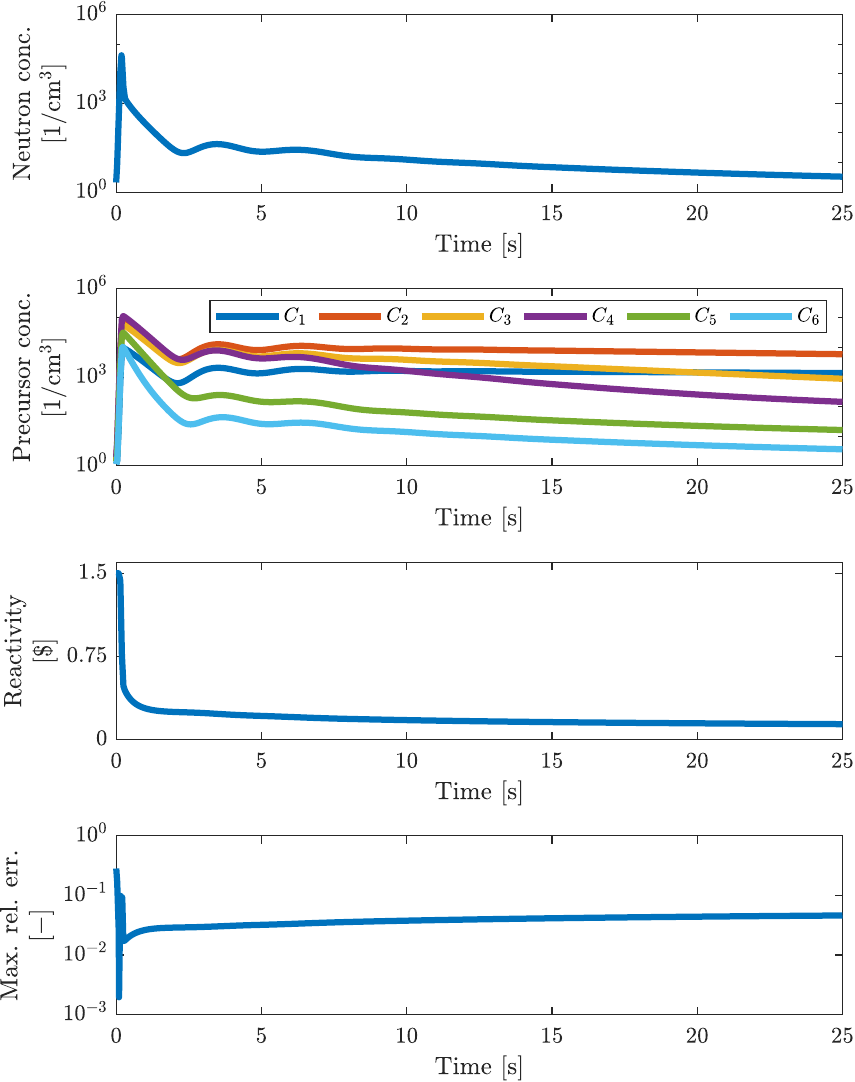}
    \caption{Estimation results for the point reactor kinetics model. Left: True and estimated kernel (top), the corresponding absolute error (second from top), the estimates of the coefficients $\{c_m\}_{m=0}^M$ (third from top), and the true (when applicable) and estimated values of $a$, $\{C_{i, 0}\}_{i=1}^n$, $\rho_0$, and $\kappa$ (bottom). Right: The concentrations (top and second from the top) and reactivity relative to $\beta$ (third from the top) for the true parameters and the maximum relative absolute difference (max. rel. err.) obtained with the estimated parameters (bottom). The colors are consistent across the left column.}
    \label{fig:nuclear:fission:nonconforming}
\end{figure*}

Fig.~\ref{fig:nuclear:fission:nonconforming} shows the estimation results for $M = 70$. The algorithm accurately estimates the kernel, the initial reactivity, and the reactivity proportionality constant, $\kappa$. However, there are some discrepancies in the estimates of the initial concentrations. Consequently, and due to the stiffness of the process, the maximum relative error of the concentrations and the reactivity (i.e., the absolute difference divided by the true value) is high initially. For most of the time interval, the largest error occurs for the reactivity. Nonetheless, the results indicate that the algorithm can also accurately estimate a bimodal kernel for a stiff multivariate set of DDEs.
\begin{table}[t]
    \tiny
    \centering
    \caption{Values of the parameters in the point reactor kinetics model.}
    \label{tab:nuclear:fission:parameters}
    \begin{tabular}{cccccc}
        \hline
        \multicolumn{6}{c}{Decay constants~[1/s]} \\
        \hline
        $\lambda_1$ & $\lambda_2$ & $\lambda_3$ & $\lambda_4$ & $\lambda_5$ & $\lambda_6$ \\
        0.0124 & 0.0305 & 0.1110 & 0.3010 & 1.1300 & 3.0000 \\
        \hline
        \multicolumn{6}{c}{Delayed neutron fractions~[--]} \\
        \hline
        $\beta_1$ & $\beta_2$ & $\beta_3$ & $\beta_4$ & $\beta_5$ & $\beta_6$ \\
        0.00021 & 0.00141 & 0.00127 & 0.00255 & 0.00074 & 0.00027 \\
        \hline
        \multicolumn{6}{c}{Other model parameters} \\
        \hline
        $\beta$~[--] & $\Lambda$~[s] & $\kappa$~[1/K] & $H$~[K\,cm$^3$/s]& $\tau_c$~[s] & $\tau_\ell$~[s] \\
        0.0065 & $5\cdot 10^{-5}$ & $5\cdot 10^{-5}$ & 0.05 & 0.5 & 3.5 \\
        \hline
        \multicolumn{6}{c}{Kernel parameters} \\
        \hline
        $\gamma_1$~[--] & $\gamma_2$~[--] & $\mu_1$~[s] & $\mu_2$~[s] & $\sigma_1$~[s] & $\sigma_2$~[s] \\
        0.6 & 0.4 & 2.5 & 5 & 0.5 & 1 \\
        \hline
    \end{tabular}
\end{table}

%% file: tex/conclusions.tex
\section{Conclusions\label{Conclusions}}
We present an algorithm for approximately identifying distributed time delays in DDEs based on discrete-time measurements.
First, we 1)~approximate the kernel using the probability density function of a mixed Erlang distribution and 2)~transform the resulting set of DDEs into ODEs using the LCT. Then, we formulate the delay identification problem as a dynamical least-squares optimization problem, where the dynamical constraints are the approximate system of ODEs. We also, simultaneously, estimate model parameters and the initial state. Finally, we use a single-shooting approach to transform the dynamical optimization problem into an NLP, and we approximate its solution using off-the-shelf software for numerical simulation and optimization.
We demonstrate, using two numerical examples, that the algorithm can identify distributed time delays with high precision. Specifically, we consider the logistic equation with a time-varying carrying capacity and a point reactor kinetics model of a molten salt nuclear fission reactor.


%% file: appendix/derivation.tex
\section{Derivation of the approximate system}\label{app:derivation}
In this appendix, we use the linear chain trick to transform the approximate system of DDEs~\eqref{eq:lct:approximate:system:DDE} to the system of ODEs~\eqref{eq:lct:approximate:system:ODE}. First, we introduce $\hat z_m: \R \rightarrow \R^{n_z}$ given by
\begin{align}\label{eq:approximation:revisited:z}
    \hat z_m(t) &= \int\limits_{-\infty}^t \hat \alpha_m(t - s) \hat r(s)\,\mathrm ds, & m &= 0, \ldots, M,
\end{align}
such that
\begin{align}\label{eq:approximation:z}
        \hat z(t) &= \int\limits_{-\infty}^t \hat \alpha^{(M)}(t - s) \hat r(s)\,\mathrm ds = \sum_{m=0}^M c_m \hat z_m(t).
\end{align}
Next, we rewrite the expressions for the time derivatives of the probability density functions of the Erlang distributions as
\begin{subequations}\label{eq:erlang:pdf:derivative}
    \begin{align}
        \dot{\hat \alpha}_0(t) &= -a b_0 e^{-at} = -a \hat \alpha_0(t), \\
        \dot{\hat \alpha}_m(t)
        &= m b_m t^{m-1} e^{-at} - a b_m t^m e^{-at} \nonumber \\
        &= a b_{m-1} t^{m-1} e^{-at} - a b_m t^m e^{-at} \nonumber \\
        &= a(\hat \alpha_{m-1}(t) - \hat \alpha_m(t)), \quad m = 1, \ldots, M,
    \end{align}
\end{subequations}
where we have used that
\begin{align}
    b_m &= \frac{a^{m+1}}{m!} = \frac{a}{m} \frac{a^m}{(m-1)!} = \frac{a}{m} b_{m-1},
\end{align}
for $m = 1, \ldots, M$. Finally, we use Leibniz' integral rule~\cite[Thm.~3, Chap.~8]{Protter:Morrey:1985} to derive the ODEs:
\begin{align}
    \dot{\hat z}_m(t)
        &= \hat \alpha_m(0) \hat r(t) + \int\limits_{-\infty}^t \dot{\hat \alpha}_m(t - s) \hat r(s)\, \mathrm ds,
\end{align}
for $m = 0, \ldots, M$. Specifically, we obtain
\begin{subequations}
    \begin{align}
        \dot{\hat z}_0(t)
        &= a \hat r(t) - a\int\limits_{-\infty}^t \hat \alpha_0(t-s) \hat r(s)\, \mathrm ds \nonumber \\
        &= a (\hat r(t) - \hat z_0(t)), \\
        \dot{\hat z}_m(t)
        &= a \bigg(\int\limits_{-\infty}^t \hat \alpha_{m-1}(t - s) \hat r(s)\,\mathrm ds \nonumber \\
        &\phantom{=}- \int\limits_{-\infty}^t \hat \alpha_m(t - s) \hat r(s)\,\mathrm ds\bigg) \nonumber \\
        &= a (\hat z_{m-1}(t) - \hat z_m(t)), \quad m = 1, \ldots, M,
    \end{align}
\end{subequations}
where we have used~\eqref{eq:erlang:pdf:derivative} and that
\begin{align}
    \hat \alpha_m(0) &=
    \begin{cases}
        a, & \text{for}~m = 0, \\
        0, & \text{for}~m = 1, \ldots, M,
    \end{cases}
\end{align}
since $b_0 = a$.
The resulting system of ODEs,
\begin{subequations}
    \begin{align}
        \dot{\hat z}_0(t) &= a (\hat r(t) - \hat z_0(t)), \\
        \dot{\hat z}_m(t) &= a (\hat z_{m-1}(t) - \hat z_m(t)), & m &= 1, \ldots, M, \\
        \hat z(t) &= \sum\limits_{m=0}^M c_m \hat z_m(t),
    \end{align}
\end{subequations}
is in the form~\eqref{eq:system:approximate:Z}--\eqref{eq:system:approximate:z}, where
\begin{align}
    \hat Z &=
    \begin{bmatrix}
        \hat z_0 \\
        \hat z_1 \\
        \vdots \\
        \hat z_M
    \end{bmatrix}.
\end{align}
\begin{rem}\label{rem:steady:state}
    If $\hat x(t) = \hat x_0$ for $t \in (-\infty, t_0]$, then
    \begin{align}\label{eq:identification:initial:conditions:zm0}
        \hat z_m(t_0)
        &= \int\limits_{-\infty}^{t_0} \hat \alpha_m(t_0 - s) \hat r(s)\,\mathrm ds \nonumber \\
        &= h(\hat x_0, p), \quad m = 0, \ldots, M.
    \end{align}
\end{rem}

%% file: appendix/numerical_simulation.tex
\section{Numerical simulation}\label{sec:numerical:simulation}
In this appendix, we describe a numerical approach for approximating the solution to initial value problems involving differential equations in the form~\eqref{eq:system:x} and distributed delays in the form~\eqref{eq:system:delay}.
First, we approximate the integral in~\eqref{eq:system:z} by
\begin{align}\label{eq:numerical:simulation:system:z}
    z(t) &= \int\limits_{-\infty}^t \alpha(t - s) r(s)\,\mathrm ds \approx \int\limits_{t - \Delta t_h}^t \alpha(t - s) r(s)\,\mathrm ds,
\end{align}
i.e., we only integrate from time $t - \Delta t_h$ to $t$. We choose $\Delta t_h \in (0, \infty)$ such that $\alpha(t) < \epsilon$ for all $t > \Delta t_h$ and some small $\epsilon \in (0, \infty)$. Furthermore, we use a fixed step size of $\Delta t \in (0, \infty)$ in the discretization of the differential equations and the integral, and we choose it such that $N_h = \Delta t_h/\Delta t \in \N$ is integer. Finally, we assume that $x(t)$ (and, consequently, also $z(t)$ and $r(t)$) are given for $t \in (-\infty, t_0]$.
Next, we discretize the differential equations using Euler's implicit method and we discretize the right-most integral in~\eqref{eq:numerical:simulation:system:z} using a right rectangle rule:
\begin{subequations}\label{eq:numerical:simulation:stiff}
    \begin{align}
        \label{eq:numerical:simulation:stiff:x}
        x_{n+1} &= x_n + f(x_{n+1}, z_{n+1}, p) \Delta t, \\
        \label{eq:numerical:simulation:stiff:z}
        z_{n+1} &= \sum_{j=0}^{N_h-1} \alpha(j \Delta t) r_{n-j+1} \Delta t, \\
        \label{eq:numerical:simulation:stiff:r}
        r_{n+1} &= h(x_{n+1}, p).
    \end{align}
\end{subequations}
For $t_n \leq t_0$, $x_n \in \R^{n_x}$ and $z_n, r_n \in \R^{n_z}$ are equal to $x(t_n)$, $z(t_n)$, and $r(t_n)$, respectively, and for $t_n > t_0$, they are approximations. Furthermore, $t_n = t_0 + n \Delta t$.
All three algebraic equations~\eqref{eq:numerical:simulation:stiff} are coupled and must be solved for $x_{n+1}$, $z_{n+1}$, and $r_{n+1}$ for given $x_n$, $z_n$, and $\{r_{n-j+1}\}_{j=1}^{N_h-1}$. We compute $x_{n+1}$ as the root of the residual function $R_n: \R^{n_x} \times \R^{n_x} \times \R^{n_p} \rightarrow \R^{n_x}$, which is given by
\begin{align}\label{eq:numerical:simulation:stiff:residual}
    R_n(x_{n+1}; x_n, p) &= x_{n+1} - x_n - f(x_{n+1}, z_{n+1}, p) \Delta t,
\end{align}
where $z_{n+1}$ and $r_{n+1}$ are treated as functions of $x_{n+1}$.
We solve the equations sequentially in a forward manner, and we use \textsc{Matlab}'s \texttt{fsolve} to approximate the roots numerically. Furthermore, we supply the analytical Jacobian of the residual function, which is
\begin{align}
    \pdiff{R_n}{x_{n+1}}(x_{n+1}; x_n, p) &= I - \Bigg(\pdiff{f}{x}(x_{n+1}, z_{n+1}, p) \nonumber \\
    &\phantom{=} + \pdiff{f}{z}(x_{n+1}, z_{n+1}, p) \pdiff{z_{n+1}}{x_{n+1}}\Bigg) \Delta t,
\end{align}
where
\begin{subequations}
    \begin{align}
        \pdiff{z_{n+1}}{x_{n+1}} &= \alpha(0) \pdiff{r_{n+1}}{x_{n+1}} \Delta t, \\
        \pdiff{r_{n+1}}{x_{n+1}} &= \pdiff{h}{x}(x_{n+1}, p),
    \end{align}
\end{subequations}
and $I \in \R^{n_x \times n_x}$ is an identity matrix.